\documentclass{ps}
\usepackage{amsmath}
\usepackage{amsfonts}
\usepackage{amssymb}
\usepackage{subfigure}
\usepackage{graphicx}

\newcommand{\blackboard}[1]{\mathbf{#1}} 
\newcommand{\Q}{\blackboard{Q}}
\newcommand{\R}{\blackboard{R}}
\renewcommand{\P}{\blackboard{P}}
\newcommand{\E}{\blackboard{E}}
\renewcommand{\L}{\blackboard{L}}

\newcommand{\calig}[1]{\mathcal{#1}}
\newcommand{\A}{\calig{A}}
\newcommand{\C}{\calig{C}}
\newcommand{\F}{\calig{F}}

\DeclareMathOperator{\sgn}{sgn}
\DeclareMathOperator{\var}{var}
\DeclareMathOperator{\supp}{supp}

\renewcommand{\d}[1]{\textup{d} #1}

\newcommand{\convps}{
    \xrightarrow[\phantom{ a.s }]{a.s}}

\newcommand{\convet}{
    \Rightarrow}

\newcommand{\delimleft}[2]{\ifcase #1\or
    \bigl#2\or %
    \Bigl#2\or %
    \biggl#2\or %
    \Biggl#2\or %
    \left#2\fi}
\newcommand{\delimright}[2]{\ifcase #1\or
    \bigr#2\or %
    \Bigr#2\or %
    \biggr#2\or %
    \Biggr#2\or %
    \right#2\fi}
\newcommand{\pa}[2][5]{
    \delimleft{#1}{(} #2 \delimright{#1}{)}}
\newcommand{\br}[2][5]{
    \delimleft{#1}{[} #2 \delimright{#1}{]}}
\newcommand{\ac}[2][5]{
    \delimleft{#1}{\{} #2 \delimright{#1}{\}}}

\newcommand{\abs}[2][5]{
    {\delimleft{#1}{\lvert} #2%
    \delimright{#1}{\rvert}}}

\newcommand{\domin}[2][5]{
    \calig{O}\! \pa[#1]{#2}}
\newcommand{\neglig}[2][5]{o\! \pa[#1]{#2}}
\newcommand{\un}{
    \boldsymbol{1}}

\newcommand{\prob}[2][5]{
    \P \br[#1]{#2}}
\newcommand{\probs}[3][5]{
    \P_{#2} \br[#1]{#3}}
\newcommand{\probc}[3][5]{
    \P \br[#1]{#2%
    \delimleft{#1}{.}%
    \vphantom{#2}\vphantom{#3}%
    \delimright{#1}{\vert} #3}}
\newcommand{\esp}[2][5]{
    \E \br[#1]{#2}}
\newcommand{\esps}[3][5]{
    \E_{#2} \br[#1]{#3}}

\newcommand{\seq}[2][0]{
    {\bigl( #2_n \bigr)}_{n \ge #1}}
\newcommand{\proc}[2][0]{
    {\bigl( #2_t \bigr)}_{t \ge #1}}

\renewcommand{\ie}{\textit{i.e. }}
\renewcommand{\le}{\leqslant}
\renewcommand{\ge}{\geqslant}
\newcommand{\ds}{\displaystyle}

\begin{document}
\title{Behavior of the Euler scheme with decreasing step in a degenerate situation}
\author{Vincent Lemaire}
\date{\today}
\address{Laboratoire d'Analyse et de Math\'{e}matiques Appliqu\'{e}es, UMR 8050,\\ Universit\'{e} de Marne-la-Vall\'{e}e, 5 boulevard Descartes, Champs-sur-Marne, \\ F-77454 Marne-la-Vall\'{e}e Cedex 2, France.} 
\begin{abstract}
The aim of this paper is to study the behavior of the weighted empirical measures of the decreasing step Euler scheme of a one-dimensional diffusion process having multiple invariant measures. This situation can occur when the drift and the diffusion coefficient are vanish simultaneously. 

As a first step, we give a brief description of the Feller's classification of the one-dimensional process. We recall the concept of attractive and repulsive boundary point and introduce the concept of strongly repulsive point. That allows us to establish a classification of the ergodic behavior of the diffusion. We conclude this section by giving necessary and sufficient conditions on the nature of boundary points in terms of Lyapunov functions.

In the second section we use this characterization to study the decreasing step Euler scheme. We give also an numerical example in higher dimension.

\end{abstract}
\maketitle

\section{Introduction and framework}
Let $I=]l,r[$ denote an open (non-trivial) interval of the real line $\R$. We consider the following stochastic differential equation 
\begin{equation} \label{dim1-eds}
	\d X_t = b(X_t) \d t + \sigma(X_t) \d B_t,
\end{equation}
where $X_0$ is a random variable taking values in $I$ and $\proc{B}$ a standard Brownian motion on $\R$. 
We assume that $b$ and $\sigma$ are continuous functions on $\bar{I}$ taking values in $\R$, and that $\sigma$ is not degenerate on $I$ \ie $\forall x \in I$, $\sigma^2(x) > 0$. Then there exists a unique solution $\proc{X}$ adapted to the completed Brownian filtration, such that $t \mapsto X_t$ is continuous on $[0, \zeta[$, where $\zeta = \inf\ac{t \ge 0, \; X_t = l \text{ or } X_t = r}$ is the explosion time of the diffusion.

In the first part, we establish a new ``ergodic classification'' for the process $\proc{X}$ in particular when $\zeta = +\infty$. More precisely, we give the behavior of the sequence of empirical measures $\proc{\nu} = \pa{\int_0^t \delta_{X_s} \d s}_{t \ge 0}$ according to the nature of the \emph{boundary points} $l$ and $r$. We characterize then the nature of the boundary points in terms of Lyapunov functions. The Lyapunov functions are usually used in high dimension, but there is a close link between these functions and the Feller's classification. This link makes it possible to more easily study the Euler scheme.  

In the second part, we study the Euler scheme with decreasing step of a diffusion $\proc{X}$ on the real line having (at least) a point $\Delta$ such that $b(\Delta) = \sigma(\Delta) = 0$. In this situation we have 
\begin{equation*}
	\forall x \in ]-\infty, \Delta[, \quad \probs{x}{X_t \in ]-\infty, \Delta]} = 1, \quad \text{and} \quad  
	\forall x \in ]\Delta, +\infty[, \quad \probs{x}{X_t \in [\Delta, +\infty[} = 1.
\end{equation*}
In fact, the process $\proc{X}$ has an ergodic behavior in $I_1 = ]-\infty, \Delta[$ or in $I_2 = ]\Delta, +\infty[$ according to the starting point $x$. But the Euler scheme is not continuous and may jump above the boundary point $\Delta$. A legitimate question is then which are the weak limit of the empirical measures of the scheme ? We answer it in some cases.

\section{Results for the time continuous process}
We first introduce the scale function and the speed measure of the diffusion process $\proc{X}$ solution of \eqref{dim1-eds}. The next two sections are adapted from classical work on the Feller classification (see for instance \cite{karlin-taylor}, \cite{karatzas-shreve}, \cite{breiman} and \cite{revuz-yor}). 

\subsection{Scale function and speed measure}
\begin{dfntn}[Scale function] A \emph{scale function} for the SDE \eqref{dim1-eds} is defined for any $c \in I$ by
\begin{equation*}
	\forall x \in I, \quad p(x) = \int_c^x \exp \pa{- \int_c^y \frac{2 b(z)}{\sigma^2(z)} \d z} \d y.
\end{equation*}
\end{dfntn}
A scale function $p$ is a strictly increasing function defined up to an affine transformation. For the sake of simplicity, we call $p$ \emph{the} scale function of the process $\proc{X}$.

We notice that the continuity of $b$ and $\sigma$, and the non-degeneracy of $\sigma$ on $I$ imply that $p$ is in $\C^2(I, \R)$. Moreover, $p$ is a solution of the following ordinary differential equation 
\begin{equation*}
	\forall x \in I, \quad b(x) p'(x) + \frac{1}{2} \sigma^2(x) p''(x) = 0 \quad \ie \quad \A p(x) = 0,
\end{equation*}
and this property characterizes it. The probability that the process starting at $x$ hits a point $a \in I$ before a point $b \in I$ is then expressed by using the scale function $p$. For any $a \in I$ we denote $T_a$ the hitting time of the one-point set $\ac{a}$ \ie $T_a = \inf \ac{t \ge 0, \; X_t = a}$, and we consider a non-trivial interval $]a, b[ \subset I$ (strictly included).
The function $u(x)$ defined on $]a, b[$ by $\forall x \in ]a, b[$, $u(x) = \probs{x}{x_{T_a \wedge T_b} = b}$ is solution to the system  
\begin{equation*}
\begin{cases}
	\A u = 0 \quad \text{on $]a, b[$} \\
	u(a) = 0 \quad \text{and} \quad u(b) = 1,
\end{cases}
\end{equation*}
so that for all $x \in ]a, b[$, 
\begin{equation} \label{def-fct-p-2}
	u(x) = \probs{x}{T_b < T_a} = \frac{p(x) - p(a)}{p(b) - p(a)}.
\end{equation}
This characterization is often used to define the scale function in a more general framework \ie for continuous strongly Markovian processes which are regular in Dynkin's sense ($\forall x \in I$, $\forall y \in I$, $\probs{x}{T_y < +\infty} > 0$) (cf. \cite{rogers-williams} or \cite{revuz-yor}). 

The following proposition gives another characterization of the scale function. 
\begin{prpstn} \label{caract-p} The process $\pa{p(x^\zeta_t)}_{t \ge 0}$ is a local martingale if and only if $p$ is the scale function.
\end{prpstn}

\begin{proof} For a proof in a more general framework, see proposition VII.3.5 in \cite{revuz-yor}. 

If $\pa{p(x^\zeta_t)}_{t \ge 0}$ is a local martingale, then for all $a < x < b$ the process $\pa{p(x^{T_a \wedge T_b}_t)}_{t \ge 0}$ is a bounded martingale and by the optional sampling theorem, we have  
\begin{equation*}
	p(x) = p(a) \probs{x}{T_a < T_b} + p(b) \probs{x}{T_b < T_a},
\end{equation*}
which implies \eqref{def-fct-p-2}.

If $p$ is the scale function, then $\A p = 0$, and by the Ito's lemma applied to $\proc{x^\zeta}$ with $p$ we deduce that $\pa{p(x^\zeta_t)}_{t \ge 0}$ is a local martingale.
\end{proof}
The above proposition is very useful because it makes it possible to consider a one-dimensional diffusion as a Brownian local martingale up to a simple transform.

Namely, the process $\proc{Y}$ defined for every $t \ge 0$ by $Y_t = p(X_t)$ satisfies the following equation
\begin{equation} \label{eds-y}
	Y_t = Y_0 + \int_0^t g(Y_s) \d B_s,
\end{equation}
where $Y_0 = p(X_0) \in p(I)$ $a.s.$ and $g$ is defined by  
\begin{equation*}
	g(y) = \begin{cases} 
	    \pa{(\sigma p') \circ p^{-1}}(y) & \text{if $y \in p(I)$}, \\
	    0 & \text{otherwise}.
	\end{cases}
\end{equation*}
The process $\proc{Y}$ can be seen as time-changed Brownian motion. The speed measure of $\proc{Y}$ evaluates how the change time affects the average time of exit from a bounded open interval. Let $\tilde{\A}$ be the generator of the process $\proc{Y}$ and $v$ the function defined on $J=]\tilde{a}, \tilde{b}[ \subset p(I)$ by $v(y) = \esps{y}{\tilde{T}_{\tilde{a}} \wedge \tilde{T}_{\tilde{b}}}$ where $\tilde{T}_z = \inf \ac{t \ge 0, \; Y_t = z}$. The function $v$ is solution to the system  
\begin{equation}
	\begin{cases} 
	\tilde{\A} v(y) = -1, \quad \tilde{a} < y < \tilde{b}, \\
	v(\tilde{a}) = v(\tilde{b}) = 0.
	\end{cases}
\end{equation}
Moreover by \eqref{eds-y} we have $\tilde{\A} v(y) = \frac{1}{2} g^2(y) v''(y)$ for every $y \in ]\tilde{a}, \tilde{b}[$ and then
\begin{equation} \label{dim1-sol-v}
	\forall y \in ]\tilde{a}, \tilde{b}[, \quad v(y) = \int_{\tilde{a}}^{\tilde{b}} G_{\tilde{a},\tilde{b}}(y, z) \frac{2}{g^2(z)} \d z,
\end{equation}
where $G_{\tilde{a}, \tilde{b}}(y, z)$ is the Green function defined by $G_{\tilde{a}, \tilde{b}}(y, z) = \frac{\pa[1]{y \wedge z - \tilde{a}} \pa[1]{\tilde{b} - y \vee z}}{\tilde{b} - \tilde{a}}$, for every $(y, z) \in ]\tilde{a}, \tilde{b}[^2$.

\begin{dfntn}[Speed measure] The speed measure of the time-changed Brownian motion $\pa{p(X_t)}_{t \ge 0}$ defined in \eqref{eds-y} is the measure $\tilde{M}$ with density $\tilde{m} = 2 g^{-2}$ with respect to the Lebesgue measure.

The speed measure of the process $\proc{X}$ is the image of $\tilde{M}$ by $p^{-1}$ and is a measure with density $m = \frac{2}{\sigma^2 p'}$ with respect to the Lebesgue measure.
\end{dfntn}

As the speed measure of $\proc{X}$ is the image of the speed measure of $\proc{Y}$ by $p^{-1}$ it follows from \eqref{dim1-sol-v} that  
\begin{equation*}
	\forall y \in ]\tilde{a}, \tilde{b}[, \quad v(y) = \esps{y}{\tilde{T}_{\tilde{a}} \wedge \tilde{T}_{\tilde{b}}} = \int_{p^{-1}(\tilde{a})}^{p^{-1}(\tilde{b})} G_{\tilde{a}, \tilde{b}}(y, p(z)) m(z) \d z.
\end{equation*}
Denoting $a = p^{-1}(\tilde{a})$ and $b = p^{-1}(\tilde{b})$, we have for every $x \in ]a, b[$ 
\begin{equation*}
	v(p(x)) = \int_a^b G_{p(a), p(b)}(p(x), p(z)) m(z) \d z.
\end{equation*}
Moreover $p$ is a one-to-one function, then it is straightforward that $\esps{p(x)}{\tilde{T}_{p(a)} \wedge \tilde{T}_{p(b)}} = \esps{x}{T_a \wedge T_b}$. Using \eqref{def-fct-p-2} we check that for all $]a, b[ \subset I$ and $x \in ]a, b[$  
\begin{equation} \label{inter-prob-mesure-vitesse}
	\esps{x}{T_a \wedge T_b} = \pa[1]{1-\probs{x}{T_b < T_a}} \int_a^x \pa[1]{p(y) - p(a)} m(y) \d y \\ + \probs{x}{T_b < T_a} \int_x^b \pa[1]{p(b) - p(y)} m(y) \d y.
\end{equation}

In addition the scale function and the speed measure provide a very useful expression for the infinitesimal generator $\A$ associated the SDE \eqref{dim1-eds}. Indeed, we easily check that 
\begin{equation} \label{expression-generateur-p-m}
	\forall f \in \C^2(I, \R), \quad \A f(x) = b(x) f'(x) + \frac{1}{2} \sigma^2(x) f''(x) = \frac{1}{m(x)} \pa{\frac{f'(x)}{p'(x)}}'.
\end{equation}

\subsection{Feller classification}
The classification of one-dimensional diffusion process is due to Feller, in particular in \cite{feller-parabolic} and \cite{feller-classif}.  In parallel, the Russian school established similar results, but the two terminologies do not always coincide.  For a comparison and a synthesis we refer to \cite{karlin-taylor}.

A first concept characterizing the behavior of the process $\proc{X}$ solution of \eqref{dim1-eds} in a neighborhood of a boundary point of $I = ]l, r[$ is the attractivity. 
\begin{dfntn}[Attractivity] A boundary point $\Delta$ ($\Delta = l$ or $\Delta = r$) is said to be \emph{attractive} if $\ds \lim_{\substack{b \rightarrow \Delta\\b \in I}} \abs{p(b)} < +\infty$. 
\end{dfntn}

The function $p$ is defined up to a strictly increasing affine transformation but the fact that the limit in $\Delta$ of $p$ is finite (or infinite) does not depend on it. Similarly, since $\abs{p(x) - p(y)} < +\infty$ for all $x, y$ and that $p$ is strictly increasing we have 
\begin{align*}
	\text{$l$ is attractive} \quad & \Leftrightarrow \quad \forall x \in I, \quad \lim_{b \rightarrow l} \pa[1]{p(x) - p(b)} < +\infty, \\
	\text{$r$ is attractive} \quad & \Leftrightarrow \quad \forall x \in I, \quad \lim_{b \rightarrow r} \pa[1]{p(b) - p (x)} < +\infty.
\end{align*}
These equivalences are sometimes used to define the attractivity.  We now show the following proposition which justifies the name of ``attractive point''.  
\begin{prpstn}\label{prop-attractif} If $\Delta$ is an attractive boundary point, then for all $a \in I$ and all $x$ in the open interval with endpoints $a$ and $\Delta$ we have 
\begin{equation*}
	\probs{x}{T_\Delta \le T_a} > 0,
\end{equation*}
with $T_\Delta$ defined by $\ds T_\Delta = \lim_{\substack{b \rightarrow \Delta \\ b \in I}} T_b$.
\end{prpstn}

\begin{proof} 
We detail the proof for the case $a < x < \Delta$. Firstly, since $\proc{X}$ is continuous we have by the mean-value theorem that the function $\pa{b \mapsto T_b}$ (for $b > x$) is strictly increasing. Thus $T_\Delta$ is the strictly increasing limit of $T_b$ when $b$ increases to~$\Delta$. Hence we have 
\begin{equation*}
	\bigcap_{\substack{b \in I \cap \Q \\ x < b < \Delta}} {}^{\downarrow}\ac{T_b < T_a} = \ac{T_\Delta \le T_a},
\end{equation*}
and $\ds \probs{x}{T_\Delta \le T_a} = \lim_{b \rightarrow \Delta} \probs{x}{T_b < T_a}$. By \eqref{def-fct-p-2} we deduce that 
\begin{equation*}
	\probs{x}{T_\Delta \le T_a} = \lim_{b \rightarrow \Delta} \frac{p(x) - p(a)}{p(b) - p(a)}.
\end{equation*}
However $\Delta$ is attractive, therefore $\ds \lim_{b \rightarrow \Delta} \pa[1]{p(b) - p(a)}< +\infty$, and the proof is complete.
\end{proof}

\begin{dfntn}[Repulsivity] A boundary point $\Delta$ is said to be \emph{repulsive} if it is not attractive, \ie $\ds \lim_{\substack{b \rightarrow \Delta \\ b \in I}} \abs{p(b)} = +\infty$.
\end{dfntn}
Since $p$ is strictly increasing and finite for any point of $I$, it is clear that
\begin{align*}
	\text{$l$ is repulsive} \quad & \Leftrightarrow \quad \forall x \in I, \quad \lim_{b \rightarrow l} \pa[1]{p(x) - p(b)} = +\infty, \\
	\text{$r$ is repulsive} \quad & \Leftrightarrow \quad \forall x \in I, \quad \lim_{b \rightarrow r} \pa[1]{p(b) - p (x)} = +\infty.
\end{align*}
We also check by proposition \ref{prop-attractif} that $\Delta$ repulsive implies that for all $a \in I$ and $x$ in the open interval of endpoints $a$ and $\Delta$, $\probs{x}{T_\Delta > T_a} = 1$.

\subsubsection{Attainability}
If $\Delta$ is an attractive boundary point, then a trajectory starting at $x \in I$ hits $\Delta$ before another point $a \in I$ with strictly positive probability. But does this event occur in a finished time?  Yes, if the point $\Delta$ is attainable.
\begin{dfntn}[Attainability] A boundary point $\Delta$ is said to be \emph{attainable} if for all $a \in  I$ and $x$ in the open interval of endpoints $a$ and $\Delta$ we have  
\begin{equation*}
	\lim_{b \rightarrow \Delta} \esps{x}{T_b \wedge T_a} < +\infty.
\end{equation*}
\end{dfntn}

An attainable boundary point is attractive, and if $\Delta$ is attractive, then $\Delta$ is attainable if and only if $\probs{x}{T_\Delta < +\infty} > 0$ (cf. lemma 6.2 in \cite{karlin-taylor}).
 
These two concepts ``attractivity'' and ``attainability'' make it possible to determine the behavior of the diffusion in an neighborhood of a boundary point.  Note that other concepts of the Feller's classification are not evoked here:  regular point (reflective, absorbent, adhesive), exit point, natural point, entrance point. 

\subsection{Behavior of empirical measures: the ergodic point of view}
Using Feller's classification, we can know the asymptotic behavior of one trajectory of the diffusion. But to establish the behavior of empirical measures and the recurrence of the process, the concept of attractivity is not precise enough. Indeed, several situations can occur: the boundary points $+\infty$ and $-\infty$ are both repulsive for the Brownian motion and for the Ornstein-Uhlenbeck process. But for the Brownian motion it is null recurrent (in dimension one) and for the O.U. process it is positive recurrent. A new concept then is introduced: ``strong repulsivity''.
\begin{dfntn} A repulsive boundary point $\Delta$ is said to be \emph{strongly repulsive} if for every $c \in I$ we have $\abs{\int_{\Delta}^c m(y) \d y} < +\infty$.
\end{dfntn}

A strongly repulsive point is a \emph{repulsive} point such that the speed measure is finite in a neighborhood of $\Delta$.  We emphasize that this concept is defined from the repulsivity. Indeed an attractive boundary point $\Delta$ may satisfy $\abs{\int_{\Delta}^c m(y) \d y} < +\infty$ for every $c \in I$ (see the following example).
\begin{xmpl} \begin{enumerate}
\item Let $I = ]0, +\infty[$ and $b$ and $\sigma$ continuous on $\bar{I}$. Furthermore assume that for every $x \in [0, 1]$, $b(x) = \frac{1}{2} \sqrt{x}$ and $\sigma(x) = c x^{3/4}$ with $c \in ]1, \sqrt{2}[$. We have 
\begin{equation*}
	\forall x \in ]0, 1[, \quad p'(x) = \exp\pa{- \int_1^x \frac{\sqrt{y}}{c^2 y^{3/2}} \d y} = x^{-\frac{1}{c^2}}, 
\end{equation*}
and since $c > 1$, for every $x \in ]0, 1[$, $p(x) = \frac{1}{1-1/c^2} x^{1-\frac{1}{c^2}} - \frac{1}{1-1/c^2}$ therefore $\ds \lim_{x \rightarrow 0} p(x) = - \frac{1}{1-1/c^2}$. The boundary point $0$ is thus attractive. 

Moreover, the speed measure is finite in a neighborhood of $0$ since
\begin{equation}
	\int_0^1 m(y) \d y = \int_0^1 \frac{2}{c^2 y^{3/2} y^{-1/c^2}} \d y = \frac{2}{c^2} \frac{1}{1/c^2 - 1/2},
\end{equation}
and $c \in ]1, \sqrt{2}[$.

\item In the case of the Ornstein-Uhlenbeck process defined on $\R$ by  
\begin{equation}
	\d X_t = -\frac{1}{2} X_t \d t + \d B_t, \quad X_0 = x \in \R,
\end{equation}
the speed measure $m(x) \d x$ is the Gaussian probability, and the boundary points $-\infty$ and $+\infty$ are thus strongly repulsive. 
\end{enumerate}
\end{xmpl}

We recall now the main ergodic result for one dimensional diffusion process $\proc{X}$. A process $\proc{X}$ is said to be \emph{recurrent} if for all $a$ and $b$ in $I$ we have $\probs{a}{T_b < +\infty} = 1$. Moreover the process is called \emph{positive recurrent} if $\esps{a}{T_b} < +\infty$ and $\esps{b}{T_a} < +\infty$, and \emph{null recurrent} otherwise.
\begin{thrm} \label{thm-ergodique-dim1} We suppose that $\proc{X}$ (solution of \eqref{dim1-eds} and with speed measure $m$) is recurrent on $I$.
If $f$ and $g$ are two non negative measurable functions such that 
\begin{equation*}
	\int_I \pa{f(x) + g(x)} m(\d x) < +\infty, \quad \int_I g(x) m(\d x) \neq 0.
\end{equation*}
Then 
\begin{equation*}
	\frac{\ds \int_0^t f(X_s) \d s}{\ds \int_0^t g(X_s) \d s} \convps \frac{\ds \int_I f(x) m(\d x)}{\ds \int_I g(x) m(\d x)}.
\end{equation*}
\end{thrm}
\begin{proof}
For a detailed proof we refer the reader to \cite{rogers-williams} or \cite{ito-mckean}. 
\end{proof}

In the sequel we will denote by $\proc{\nu}$ the empirical measures of the diffusion, \ie 
\begin{equation*}
	\forall t > 0, \quad \nu_t(\d x) = \frac{1}{t} \int_0^t \delta_{X_s}(\d x) \d s.
\end{equation*}

By the above theorem and the concept of strongly attractive boundary point we establish the following classification of the ergodic behavior of the diffusion.
\begin{thrm} \label{thm-classif} We recall that $\zeta = \inf \ac{t \ge 0, X_t = l \text{ or } X_t = r}$. Then 
\begin{itemize}
    \item if $l$ is attractive and $r$ is repulsive then $x^\zeta_t \convps l$,
    \item if $l$ and $r$ are attractive then
    \begin{equation*}
	\prob{\lim_{t \rightarrow +\infty} x^\zeta_t = l} = 1 - \prob{\lim_{t \rightarrow +\infty} x^\zeta_t = r} = \frac{p(r^-) - p(X_0)}{p(r^-) - p(l^+)}.
    \end{equation*}
    \item if $l$ and $r$ are repulsive then the diffusion is recurrent and does not explode ($\zeta = +\infty$ $a.s.$). More precisely 
	\subitem -- if $l$ and $r$ are strongly repulsive (\ie the speed measure is finite) then the diffusion is positive recurrent and $\nu_t \convet \nu$ $a.s.$ where $\nu$ is the normalized speed measure. Moreover    
	\begin{equation*}
		\forall f \in \L^1(\nu), \quad \frac{1}{t} \int_0^t f(X_s) \d s \convet \int_{\R} f(x) \nu(\d x) \quad a.s. 
	\end{equation*}
	\subitem -- if $l$ is strongly repulsive and $r$ is (simply) repulsive then the diffusion is null recurrent and if the empirical measures are tight we have  
	\begin{equation*}
	    \frac{1}{t} \int_0^t \delta_{X_s} \d s \convet \delta_r \quad a.s.
	\end{equation*}
	\subitem -- if $l$ and $r$ are not strongly repulsive then the diffusion is null recurrent, and if the empirical measures are tight then any weak limit of $\pa{\frac{1}{t} \int_0^t \delta_{X_s} \d s}_{t \ge 1}$ is a measure with support $\{l, r\}$.
\end{itemize}
\end{thrm}

We first prove the following lemma.
\begin{lmm} \label{dim1-lem-rec-pos} If the two boundary points $l$ and $r$ are repulsive then the diffusion is positive recurrent if and only if its speed measure is finite. 
\end{lmm}
\begin{proof} By definition, the diffusion is positive recurrent if and only if for all $a$ and $b$ in $I$, $\esps{a}{T_b} < +\infty$ and $\esps{b}{T_a} < +\infty$. Let $l < a < b < r$. By symmetry it is sufficient to prove that 
\begin{equation} \label{lem-recpos-a-montrer}
	\esps{a}{T_b} < +\infty \quad \Leftrightarrow \quad \int_l^a m(y) \d y < +\infty.
\end{equation}

Firstly, $l$ is repulsive and $T_l = \lim_{x \rightarrow l} T_x$ therefore $\esps{a}{T_b} = \esps{a}{T_b \wedge T_l} = \lim_{x \rightarrow l} \esps{a}{T_b \wedge T_x}$. 
Moreover by \eqref{inter-prob-mesure-vitesse} we have $\forall x \in ]l, a[$,   
\begin{equation*}
	\esps{a}{T_b \wedge T_x} = \probs{a}{T_b < T_x} \int_a^b \pa{p(b) - p(y)} m(y) \d y + \probs{a}{T_x \le T_b} \int_x^a \pa{p(y) - p(x)} m(y) \d y,
\end{equation*}
and since $\int_a^b \pa{p(b) - p(y)} m(y) \d y$ is finite and does not depend on $x$, the limit when $x$ tends to $l$ of $\esps{a}{T_b \wedge T_x}$ is finite if and only if  
\begin{equation*}
	\lim_{x \rightarrow l} \pa{\probs{a}{T_x \le T_b} \int_x^a (p(y) - p(x)) m(y) \d y} < +\infty.
\end{equation*}
As $\probs{a}{T_x \le T_b} = \frac{p(b)-p(a)}{p(b)-p(x)}$ we have 
\begin{align*}
	\probs{a}{T_x \le T_b} \int_x^a (p(y)-p(x)) m(y) \d y & = \pa{p(b)-p(a)} \int_x^a \frac{p(y)-p(x)}{p(b)-p(x)} m(y) \d y, \\
	& = \pa{p(b)-p(a)} \int_x^a \probs{y}{T_b < T_x} m(y) \d y,
\end{align*}
and it follows that $\esps{a}{T_b} < +\infty$ if and only if $\ds \lim_{x \rightarrow l} \int_x^a \probs{y}{T_b < T_x} m(y) \d y < +\infty$. 
Since $T_x$ strictly increases to $T_l$, $\probs{y}{T_b < T_x}$ increases to $\probs{y}{T_b < T_l} = 1$ because $l$ is repulsive. The monotone convergence theorem yields \eqref{lem-recpos-a-montrer}.
\end{proof}

\begin{proof}[Proof of Theorem \ref{thm-classif}] The first two items are proved in \cite{karatzas-shreve} (Proposition 5.22). We recall the proof of the first item.
\begin{itemize}
\item Suppose that $l$ is attractive and $r$ is repulsive. By definition of the scale function we have for all $l < a < x < b < r$ 
\begin{equation*}
	\probs{x}{\inf_{0 \le t < \zeta} X_t \le a} \ge \probs{x}{x_{T_a \wedge T_b} = a} = \frac{p(b) - p(x)}{p(b) - p(a)}.
\end{equation*}
Increasing $b$ to $r$ we obtain $\probs{x}{\inf_{0 \le t < \zeta} X_t \le a} = 1$ since $r$ is repulsive.
The limit when $a$ decreases to $l$ also gives 
\begin{equation} \label{dim1-preuve-point-repulsif}
	\probs{x}{\inf_{0 \le t < \zeta} X_t = l}  = 1.
\end{equation}
On the other hand  
\begin{equation*}
	\probs{x}{\sup_{0 \le t < \zeta} X_t = b} \le \probs{x}{x_{T_l \wedge T_b} = b} = \lim_{a \rightarrow l} \frac{p(x) - p(a)}{p(b) - p(a)},  
\end{equation*}
and taking the limit when $b$ increases to $r$ we obtain $\sup_{0 \le t < \zeta} X_t < r$ $a.s.$ It remains to prove that the process $\proc{x^\zeta}$ is almost surely convergent. Since $\pa{p(x^\zeta_t)}_{t \ge 0}$ is a local martingale and $l$ is an attractive boundary point, the process $\pa{p(x^\zeta_t) - \lim_{a \rightarrow l} p(a)}_{t \ge 0}$ is a positive continuous local martingale. By Fatou's lemma it is a positive continuous super-martingale which is then almost surely convergent. We conclude using the continuity of $p^{-1}$.

\item If $l$ and $r$ are repulsive then in the same way that we obtain \eqref{dim1-preuve-point-repulsif} we have 
\begin{equation*}
	\inf_{0 \le t < \zeta} X_t = l \quad a.s. \quad \text{and} \quad \sup_{0 \le t < \zeta} X_t = r \quad a.s.
\end{equation*}
The diffusion is thus recurrent on $I$ and $\zeta = +\infty$ $a.s.$ Moreover, by Lemma \ref{dim1-lem-rec-pos} we know that the recurrence is positive if and only if the speed measure is finite.

\subitem -- If the two boundary points are repulsive then the speed measure is finite and by Theorem \ref{thm-ergodique-dim1} we have 
\begin{equation*}
        \forall f \in \L^1(\nu), \quad \frac{1}{t} \int_0^t f(X_s) \d s \convet \int_{\R} f \d \nu \quad a.s.
\end{equation*}
where $\nu$ is the normalized speed measure.

\subitem -- If $l$ is strongly repulsive and $r$ is repulsive then the diffusion is null recurrent. Considering an increasing sequence of continuous functions with compact support $\pa{g_n(x)}_{n \ge 1}$ such that $g_n(x) \rightarrow 1$ and $\forall n \ge 1$, $\int_\R g_n(x) m(\d x) \neq 0$, we obtain by Theorem \ref{thm-ergodique-dim1} 
\begin{equation} \label{dim1-rec-nulle}
	\forall f \in \L^1(m), \quad \frac{1}{t} \int_0^t f(X_s) \d s \convps 0.
\end{equation}
On the other hand, we consider a sub-sequence $\pa{\nu_{a(t)}}_{t \ge 0}$ of $\proc{\nu}$ converging to a measure $\nu$ (the empirical measures are tight). Let $f$ be a continuous function with compact support such that $\supp(f) \subset [l, r[$. As $\nu_{a(t)} \convet \nu$ we have  
\begin{equation*}
	\frac{1}{a(t)} \int_0^{a(t)} f(X_s) \d s \convps \int f \d \nu.
\end{equation*}
The boundary point $l$ is strongly repulsive and $\supp(f) \subset [l, r[$, thus $f$ is integrable with respect to $m$. Hence \eqref{dim1-rec-nulle} implies $\int f \d \nu = 0$. The interval $[l, r[$ satisfies 
\begin{equation*}
	\forall f \in \C_c(\bar{I}), \; \supp(f) \subset [l r[, \Rightarrow \nu(f) = 0,
\end{equation*}
therefore $\supp(\nu) = \ac{r}$. Since $\nu$ is normalized we have $\nu = \delta_r$. 

\subitem -- In the same way, if the two boundary points are strongly repulsive then any weak limit of $\proc{\nu}$ is a measure with support $\{l, r\}$.
\end{itemize}
\end{proof}

\subsection{Attractivity and Lyapunov function}
In order to study the behavior of the Euler scheme (with decreasing step) we establish a link between the concepts of attractivity, repulsivity and strong repulsivity, and the Lyapunov functions. Let us notice that there is few work which relates to this subject. Indeed, the Lyapunov functions are useful in high dimension and the Feller's classification is established for one-dimensional processes.

In the sequel, we will denote by $J_\Delta \subset I$ an open (non-trivial) interval included in $I$ with endpoint $\Delta$.

\begin{prpstn} \label{lien-lyapounov-feller} Let $\Delta$ be a boundary point (finite or infinite, left endpoint or right endpoint) of $I$. The following statements are equivalents  
\begin{enumerate}
\item $\Delta$ is a repulsive boundary point of $I$ if and only if there exists a neighborhood $J_\Delta \subset I$ of $\Delta$ and a strictly monotone function $V \in \C^2(J_\Delta, \R_+)$ such that 
\begin{equation*}
	\lim_{x \rightarrow \Delta} V(x) = +\infty \quad \text{and} \quad \forall x \in J_\Delta, \quad \A V(x) \le 0.
\end{equation*}
\item $\Delta$ is a strongly repulsive boundary point of $I$ if and only if there exists a neighborhood $J_\Delta \subset I$ of $\Delta$ and a strictly monotone function $V \in \C^2(J_\Delta, \R_+)$ such that   
\begin{equation*}
	\exists \varepsilon > 0, \; \forall x \in J_\Delta, \quad \A V(x) \le -\varepsilon. 
\end{equation*}
\item $\Delta$ is a attractive point of $I$ if and only if there exists a neighborhood $J_\Delta \subset I$ of $\Delta$ and a strictly monotone function $V \in \C^2(J_\Delta, \R_+)$ such that   
\begin{equation*}
	\sup_{x \in J_\Delta} V(x) = V(\Delta) < +\infty \quad \text{and} \quad \forall x \in J_\Delta, \quad \A V(x) \ge 0. 
\end{equation*}
\end{enumerate}
\end{prpstn}

\begin{proof}
We give the proof when $\Delta$ is a right endpoint of $I$. Then $J_\Delta$ is an interval $]c, \Delta[$ with $c \in I$.
\begin{itemize}
\item -- We suppose that there exists a neighborhood $J_\Delta$ of $\Delta$ and a function $V \in \C^2(J_\Delta, \R_+)$ such that $\lim_{x \rightarrow \Delta} V(x) = +\infty$ and $\A V \le 0$ on $J_\Delta$. For every $x \in J_\Delta$ we have  
\begin{equation*}
	\A V(x) = \frac{1}{m(x)} \pa{\frac{V'(x)}{p'(x)}}' \le 0
\end{equation*}
hence $V'/p'$ is decreasing on $J_\Delta$. There also exists $C > 0$ such that $V'(x) \le C p'(x)$ for every $x \in ]c, \Delta[$ since $p'>0$. It follows that $\ds \lim_{x \rightarrow \Delta} p(x) - p(c) = +\infty$ because $V$ tends to infinity in $\Delta$.

\noindent -- Conversely we must find the good Lyapunov function $V$. Let $c > 0$ be such that $p(c) > 0$ ($c$ exists because $\Delta$ is repulsive). Since $p$ is strictly increasing we have for every $x \in ]c, \Delta[$, $p(x) > p(c) > 0$. We also define the function $V$ on $]c, \Delta[$ by 
\begin{equation*}
	\forall x \in ]c, \Delta[, \quad V(x) = p(x) - p(c).
\end{equation*}
The point $\Delta$ is repulsive thus $V$ increases to infinity when $x$ tends to $\Delta$. Moreover $V \in \C^2(]c, \Delta[, \R_+)$ and $\A V = 0$.

\item -- Let $V \in \C^2(J_\Delta, \R_+)$ be such that $\lim_{x \rightarrow \Delta} V(x) = +\infty$ and $\varepsilon > 0$ such that $\A V \le - \varepsilon$. Thus we have  
\begin{equation} \label{base AV < -1}
	\int_{J_\Delta} \A V(x) m(x) \d x \le - \varepsilon \int_{J_\Delta} m(x) \d x,
\end{equation}
and since $\ds \A V(x) = \frac{1}{m(x)} \pa{\frac{V'(x)}{p'(x)}}'$ we obtain for every $c \in {J_\Delta}$,  
\begin{equation} \label{eg AV < -1}
	\int_c^\Delta \A V(x) m(x) \d x = \lim_{x \rightarrow \Delta} \frac{V'(x)}{p'(x)} - \frac{V'(c)}{p'(c)}.
\end{equation}
By \eqref{base AV < -1} we derive that  
\begin{equation*}
	\int_c^\Delta m(x) \d x \le \frac{1}{\varepsilon}\pa{\frac{V'(c)}{p'(c)} - \lim_{x \rightarrow \Delta} \frac{V'(x)}{p'(x)}}.
\end{equation*}
As the functions $V$ and $p$ are increasing on $J_\Delta$ we have $\ds \lim_{x \rightarrow \Delta} \frac{V'(x)}{p'(x)} \ge 0$, which gives $\int_c^\Delta m(x) \d x \le C$ (\ie $\Delta$ strongly repulsive).

\noindent -- Conversely we assume that $\Delta$ is strongly repulsive. Let $c \in I$ and $V$ the function defined on $]c, \Delta[$ by
\begin{equation*}
	\forall x \in ]c, \Delta[, \quad V(x) = \int_c^x \pa{p'(y) \int_y^\Delta m(z) \d z} \d y.
\end{equation*}
It is clear that $V \in \C^2(]c, \Delta[, \R_+)$ and that for every $x \in ]c, \Delta[$, $V'(x) = p'(x) \int_x^\Delta m(z) \d z$. Moreover   
\begin{equation*}
	\forall x \in ]x, \Delta[, \quad \A V(x) = -1
\end{equation*}

\item We prove (3) in the same manner as (1). For the converse we consider the function $V(x) = p(x)-p(c)$ on $]c, \Delta[$ with $c$ such that $p(c) > 0$.
\end{itemize}
\end{proof}

The above Proposition provides a useful criterion to know the nature of a boundary point. However, when the boundary point $\Delta$ is finite, a ``natural'' Lyapunov function has a minimum at $\Delta$ and it is not the case of $V$. But an easy transform allows us to obtain this property. This is the interest of the following Corollary. 

\begin{crllr} \label{dim1-cor-caract} Let $\Delta$ a boundary point of $I$. 
\begin{enumerate}
\item $\Delta$ is a repulsive boundary point of $I$ if and only if there exists a neighborhood $J_\Delta \subset I$ of $\Delta$ and a strictly monotone function $v \in \C^2(\bar{J}_\Delta, \R_+)$ satisfying $v(\Delta) = 0$, such that   
\begin{equation} \label{condition-repulsif}
	\forall x \in J_\Delta, \quad \A v (x) \ge \frac{1}{2} \sigma^2(x) \frac{(v'(x))^2}{v(x)}, 
\end{equation}
\item $\Delta$ is a strongly repulsive boundary point of $I$ if and only if there exists a neighborhood $J_\Delta \subset I$ of $\Delta$ and a strictly monotone function $v \in \C^2(\bar{J}_\Delta, \R_+)$ having a minimum at $\Delta$, such that  
\begin{equation} \label{condition-fort-repulsif}
	\exists \varepsilon > 0, \quad \forall x \in J_\Delta, \quad \A v (x) \ge \frac{1}{2} \sigma^2(x) \frac{(v'(x))^2}{v(x)} + \varepsilon v(x),
\end{equation}
\item $\Delta$ is an attractive boundary point of $I$ if and only if there exists a neighborhood $J_\Delta \subset I$ of $\Delta$ and a strictly monotone function $v \in \C^2(\bar{J}_\Delta, \R_+)$ having a minimum at $\Delta$, such that 
\begin{equation} \label{condition-attractif}
	\forall x \in J_\Delta, \quad \A v (x) < \frac{1}{2} \sigma^2(x) \frac{(v'(x))^2}{v(x)}. 
\end{equation}
\end{enumerate}
\end{crllr}

\begin{proof} Let $J_\Delta$ a neighborhood of $\Delta$ strictly included in $I$. We consider the case in which $\Delta$ is the right endpoint of $I$ \ie $J_\Delta = ]c, \Delta[$ with $c \in I$. 
We define first for $L \ge 0$ the function $\phi_{L}$ by  
\begin{align*}
	\phi_{L} :[0, \exp(L)[ & \rightarrow \R_+, \\ 
	x & \mapsto - \ln(x) + L.
\end{align*}
The function $\phi_{L}$ is a strictly decreasing one-to-one $\C^{\infty}$ function. Moreover, for every $\C^2$ function $v$ with values in $[0, \exp(L)[$ we have 
\begin{equation} \label{corollaire-critere-relation-fond}
	\A \pa{\phi_{L} \circ v} = - \frac{\A v}{v} + \frac{1}{2} \sigma^2 \frac{(v')^2}{v^2}.
\end{equation}
\begin{itemize}
\item -- We assume that there exists $v \in \C^2([c,\Delta], \R_+)$ strictly monotone satisfying $v(\Delta) = 0$ and \eqref{condition-repulsif}. We also define on $]c, \Delta[$ the function $V$ by  $\forall x \in ]c, \Delta[$, $V(x) = \pa{\phi_{L} \circ v}(x)$ with $L = \ln(v(c))$. It is a strictly monotone function which tends to infinity in $\Delta$. By \eqref{condition-repulsif} and \eqref{corollaire-critere-relation-fond} we obtain $\A V \le 0$ on $]c, \Delta[$. The proposition \eqref{lien-lyapounov-feller} implies that $\Delta$ is repulsive.

\noindent -- Conversely if $\Delta$ is repulsive then there exists $V \in \C^2(]c, \Delta[, \R_+)$ strictly monotone which goes to $+\infty$ when $x$ tends to $\Delta$. We also define $v = \phi_{L}^{-1} \circ V$ on $]c, \Delta[$ with $L = \inf_{x \in ]c, \Delta[} V(x)$, and we extend if by continuity on $[c, \Delta]$ letting $v(c) = 1$ and $v(\Delta) = 0$. By $\A V \le 0$ and \eqref{corollaire-critere-relation-fond} we have $\A v \ge \frac{1}{2} \sigma^2 \frac{(v')^2}{v}$ on $]c, \Delta[$.  

\item -- Let $v$ a strictly monotone function on $[c, \Delta]$ such that $v(\Delta) = 0$. For $\theta \ge 0$ and $L = \ln(v(c) + \theta)$ we define $V(x) = (\phi_{L, \theta} \circ v)(x)$ for every $x \in ]c, \Delta[$. This function $V$ is strictly monotone on $]c, \Delta[$ and admits a limit (finite or not) when $x$ tends to $\Delta$. From \eqref{condition-fort-repulsif} and \eqref{corollaire-critere-relation-fond} we deduce that 
\begin{equation*}
	\forall x \in ]c, \Delta[, \quad \A V(x) \le - \varepsilon.
\end{equation*}
\noindent -- Conversely we consider the function $v = \phi^{-1}_L \circ V$ on $]c, \Delta[$ with $L = \inf_{x \in ]c, \Delta[} V(x)$ and we extend if by continuity letting $v(c) = 1$ and $v(\Delta) = \lim_{x \rightarrow \Delta} \exp(- V(x))$. 

\item The proof is similar to the two firsts items.
\end{itemize}
\end{proof}

\begin{rmrk} \label{lien_ODE} If $\Delta$ is such that $b(\Delta) = \sigma(\Delta) = 0$, then $\Delta$ is a critical point for the equation $u'=b(u)$. It is worth noting that there exists a link between the nature of the critical point $\Delta$ for the ODE $u'=b(u)$ and the nature of the boundary point $\Delta$ for the SDE.  

Indeed, if $\Delta$ is a stable critical point then there exists a Lyapunov function $F \in \C^2$ such that $F'b(u) < 0$ for every $u$ in a neighborhood of $\Delta$. If $F'/F$ is decreasing, the above corollary implies that $\Delta$ is an attractive boundary point for the SDE.

If $\Delta$ is an unstable point for the ODE $u'=b(u)$, it may be repulsive, strongly repulsive or attractive for the SDE, as shown in the following example.
\end{rmrk}

\begin{xmpl} \label{exemple-}
We consider, like in \cite{pages-esaim}, the function $V:\R \rightarrow \R_+$ defined by   
\begin{equation*}
	V(x) = \begin{cases} 
	    \pa[1]{x - 3 \sgn(x)}^2 & \text{if $\abs{x} \ge 3$}, \\ 
	    \frac{1}{72} (x^2 - 9)^2 & \text{if $\abs{x} \le 3$}
	    \end{cases}
	    \quad \text{and} \quad 
	b(x) = \begin{cases} 
	    - 2 \pa[1]{x - 3 \sgn(x)}& \text{if $\abs{x} \ge 3$}, \\ 
	    - \frac{1}{18} x^3 + \frac{1}{2} x & \text{if $\abs{x} \le 3$}
	    \end{cases},
\end{equation*}
and $b = -V'$. The ordinary differential equation $u' = b(u)$ has 3 critical points: $-3$, 0 and $3$. The points $-3$ and $3$ are stable and the point $0$ is unstable.  
\begin{figure}[ht!]
    \begin{minipage}[c]{0.46\linewidth}
	\input{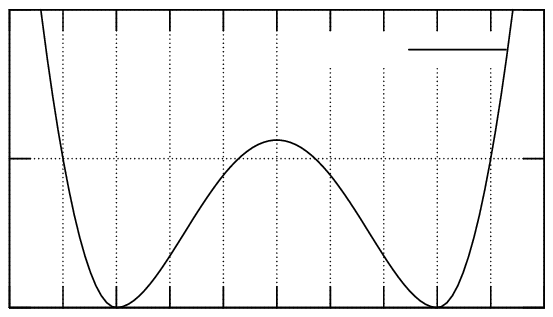}
	\caption{Fonction $V$}
    \end{minipage}
    \hfill
    \begin{minipage}[c]{0.46\linewidth}
	\input{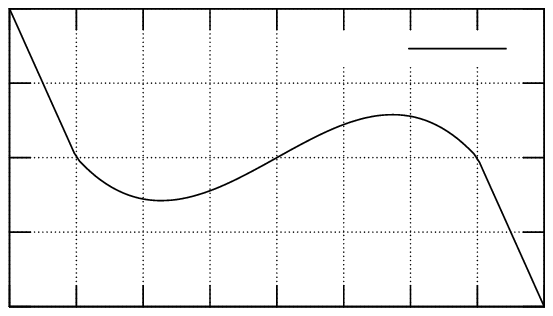}
	\caption{Drift $b = -V'$}
    \end{minipage}
\end{figure}

Let $c \in ]0, 2[$ a parameter and $\sigma$ defined by $\sigma(x) = c x$. We consider the process $\proc{X}$ solution of the SDE $\d X_t = b(X_t) \d t + \sigma(X_t) \d B_t$. It is clear that the point $0$ is a boundary point for $\proc{X}$. Moreover we check that 
\begin{equation*}
	\A V(x) = - \pa{4 - c^2} x^2 + 12 \sgn(x) \quad \text{if $\abs{x} \ge 3$}.
\end{equation*}
By Proposition \ref{lien-lyapounov-feller} the points $-\infty$ and $+\infty$ are then strongly repulsive. 

On the other hand, we use Corollary \ref{dim1-cor-caract} to determine the nature of the boundary point $0$ according to $c$. Assume that $I = ]0, + \infty[$ and let $v$ the function defined on $[0, +\infty[$ by $v(x) = x^2$. We have 
\begin{equation}
	\forall x \in ]0, 3[, \quad \A v(x) = (1+c^2) x^2 - \frac{1}{9} x^4 \quad \text{and} \quad \frac{1}{2} \sigma^2(x) \frac{(v'(x))^2}{v(x)} = 2 c^2 x^2, 
\end{equation}
and then  
\begin{itemize}
\item if $c < 1$ the boundary point $0$ is \emph{strongly repulsive} (for $]0, +\infty[$ and $]-\infty, 0[$ by symmetry). Indeed the condition \eqref{condition-fort-repulsif} is satisfied with $\varepsilon = \frac{1-c^2}{2}$ and $J_\Delta = \bigl]0, 3 \sqrt{\frac{1-c^2}{2}} \bigr[$.  
\item if $c > 1$ it is easy to check that the boundary point $0$ is \emph{attractive}. 
\item if $c = 1$ we consider the function $v(x) = x \exp(x)$ and we check that $0$ is a \emph{repulsive} boundary point. 
\end{itemize}
Thus the nature of the boundary point $0$ (which is always stable for the ODE $u'=b(u)$) may change according to $c$. By Theorem \ref{thm-classif} the ergodic behavior of $\proc{X}$ is the following:
\begin{itemize}
\item if $c \in ]1, 2[$ then $X_t \convps 0$ (for every starting point $X_0$),
\item if $c = 1$ then $\frac{1}{t} \int_0^t \delta_{X_s}(\d x) \convet \delta_0$,
\item if $c < 1$ then for every $f \in \L^1(m)$ 
\begin{equation*}
	\frac{1}{t} \int_0^t f(X_s) \d s \convps
	\begin{cases} 
		\int f d \nu_- & \text{if $X_0 \in ]0, +\infty[$} \\
		f(0)  & \text{if $X_0 = 0$} \\
		\int f d \nu_+ & \text{if $X_0 \in ]-\infty, 0[$}
	\end{cases}
\end{equation*}
with $\nu_-$ the invariant measure of $\proc{X}$ on $]-\infty, 0[$ and $\nu_+$ the invariant measure on $]0, +\infty[$.
\end{itemize}
\end{xmpl}

\section{Behavior of the Euler scheme with decreasing step}
We consider now the Euler scheme $\seq{X}$ built using a positive sequence $\seq[1]{\gamma}$ going to 0. We assume that $\seq[1]{\gamma}$ satisfies $\lim_n \sum_{k = 1}^n \gamma_k = +\infty$ and we denote $\Gamma_n = \sum_{k=1}^n \gamma_k$. The inhomogeneous Markov chain $\seq{X}$ is defined by 
\begin{equation*}
	X_{n+1} = X_n + \gamma_{n+1} b(X_n) + \sqrt{\gamma_{n+1}} \sigma(X_n) U_{n+1},
\end{equation*}
with $\seq[1]{U}$ a real white noise \ie a sequence of \emph{i.i.d.} random variables such that $\esp{U_1} = 0$ and $\var(U_1) = 1$. Furthermore, we assume that $U_1$ is a generalized Gaussian (cf. \cite{stout}) \ie such that 
\begin{equation*}
	\exists \kappa > 0, \; \forall \theta \in \R, \quad \esp[1]{\exp(\theta U_1)} \le \exp \pa{\frac{\kappa \abs{\theta}^2}{2}}.
\end{equation*}
For example $U_1$ is a standard Gaussian or a Bernoulli random variable. A consequence of the generalized Gaussian property is the following 
\begin{equation} \label{queue-sous-gauss}
	\forall a \ge 0, \quad \prob[1]{\abs{U_1} \ge a} \le \exp \pa{- \frac{a^2}{2 \kappa}}.
\end{equation}

We consider a diffusion $\proc{X}$ on the real line with $b$ and $\sigma$ continuous on $\R$. Moreover we assume that $b$ and $\sigma$ have sublinear growth \ie 
\begin{equation} \label{dim1-conditions-b-sigma}
	\exists C_b > 0, \quad \abs{b}^2 \le C_b (1+\abs{x}), \quad \text{and} \quad \exists C_\sigma > 0, \quad \abs{\sigma}^2 \le C_\sigma (1+\abs{x}). 
\end{equation}
In the sequel $\Delta$ denotes a finite point of $\R$ such that $b(\Delta) = \sigma(\Delta) = 0$, and $J_\Delta$ denotes an open interval with endpoint $\Delta$ ($J_\Delta = ]\Delta, \Delta+\varepsilon[$ if $\Delta$ is a left endpoint of $I = ]\Delta, +\infty[$ and $J_\Delta = ]\Delta-\varepsilon, \Delta[$ if it is a right endpoint of $I = ]-\infty, \Delta[$).

\subsection{Euler scheme}

We prove the following Theorem which gives the behavior of one trajectory of the Euler scheme $\seq{X}$ in this degenerate situation.  

\begin{thrm}\label{dim1-thm-euler} We assume that $\sigma$ satisfies $\sigma(x) \neq 0$ for every $x \in ]-\infty, \Delta[ \cup ]\Delta, +\infty[$ and that there exists $\mathcal{U}_\Delta = ]\Delta-\varepsilon, \Delta+\varepsilon[$ with $\varepsilon > 0$, and a convex function $v \in \C^2(\mathcal{U}_\Delta, \R_+)$ satisfying $v(\Delta) = 0$ and
\begin{equation} \label{hypotheses-v}
	\forall x \in \mathcal{U}_\Delta, \quad (v' b)(x) \ge 0 \quad \text{and} \quad \exists c_\sigma > 0, \forall x \in \mathcal{U}_\Delta, \quad \abs{(v'\sigma)(x)} \le c_{\sigma} v(x).
\end{equation}
If the step sequence $\seq[1]{\gamma}$ satisfies $\forall C > 0$, $\sum_{n \ge 1} \exp \pa{ -\frac{C}{\gamma_n} } < +\infty$ then the Euler scheme jump above $\Delta$ a finite number of times \ie 
\begin{equation*}
	\prob[1]{\exists n_0 \ge 0, \; \forall n \ge n_0, \; X_n \in ]-\infty, \Delta[} + \prob[1]{\exists n_0 \ge 0, \; \forall n \ge n_0, \; X_n \in ]\Delta, +\infty[} = 1.
\end{equation*}
\end{thrm}

We first prove the following lemma.
\begin{lmm} \label{lemme-prelim} We assume there exists a convex function $v \in \C^2(\bar{J}_\Delta, \R_+)$ satisfying $v(\Delta) = 0$ and  
\begin{equation} \label{hypotheses-v-lem}
	\forall x \in \bar{J}_\Delta, \quad (v'b)(x) \ge 0 \quad \text{and} \quad \exists c_\sigma > 0, \; \forall x \in J_\Delta, \quad \abs{(v'\sigma)(x)} \le c_{\sigma} v(x).
\end{equation}
Then, on the event $\ac{X_n \in \bar{J}_{\Delta}}$ 
\begin{equation*}
	\probc[1]{\Delta \in (X_n, X_{n+1})}{\F_n} \le \exp \pa{-\frac{1}{c_\sigma^2 \gamma_{n+1}}}.
\end{equation*}
\end{lmm}

\begin{proof}
We assume that $\Delta$ is a left endpoint of $I$ and we denote by $A_{n+1}$ the event $\ac{\Delta \in (X_n, X_{n+1})}$ (the geometric segment with endpoints $X_n$ and $X_{n+1}$). Since $v$ is continuous and $v(\Delta) = 0$ it is clear that 
\begin{align} 
	A_{n+1} &= \ac{ \exists t \in [0, 1], \; X_n + t (X_{n+1} - X_n)) = \Delta}, \notag \\
	& \subset \ac{ \exists t \in [0, 1], \; v(X_n + t (X_{n+1} - X_n)) = 0}.
	\label{algo_inclusion_An}
\end{align}
Consider $t \in [0,1]$ such that $\Delta = X_n + t(X_{n+1} - X_n)$. As $v$ is $\C^2$ on $\bar{J}_\Delta$, Taylor's formula gives  
\begin{equation*}
	v(X_n + t (X_{n+1}-X_n)) = v(X_n) + v'(X_n) t (X_{n+1} - X_n) + \frac{v''(\xi_{n+1})}{2} t^2 (X_{n+1} - X_n)^2,
\end{equation*}
with $\xi_{n+1} \in ]\Delta, X_n[$. The convexity of $v$ implies 
\begin{equation*}
	0 = v(X_n + t (X_{n+1}-X_n)) \ge v(X_n) + t \gamma_{n+1} (v' b)(X_n) + t \sqrt{\gamma_{n+1}} (v'\sigma)(X_n) U_{n+1}.
\end{equation*}
Since $v'b \ge 0$ on $\bar{J}_\Delta$ we have from \eqref{algo_inclusion_An} 
\begin{align*}
	A_{n+1} \cap \ac{X_n \in \bar{J}_\Delta} &\subset \ac{\exists t \in [0,1], \; v(X_n) + t \sqrt{\gamma_{n+1}} (v'\sigma)(X_n) U_{n+1} \le 0} \cap \ac{X_n \in \bar{J}_\Delta}, \\
	& \subset \ac{\exists t \in [0,1], \; t \sqrt{\gamma_{n+1}} \abs{(v'\sigma)(X_n) U_{n+1}} \ge v(X_n)} \cap \ac{X_n \in \bar{J}_\Delta}.
\end{align*}
Hence for any $n \ge 0$, we have on the event $\ac{X_n \in \bar{J}_\Delta}$   
\begin{align*}
	\probc{A_{n+1}}{\F_n} &\le \probc{\abs{U_{n+1}} \ge \frac{v(X_n)}{\sqrt{\gamma_{n+1}} \abs{(v' \sigma)(X_n)}}}{\F_n},  \\
	& \le \probc{\abs{U_{n+1}} \ge \frac{1}{c_\sigma \sqrt{\gamma_{n+1}}}}{\F_n},
\end{align*}
by the domination assumption \eqref{hypotheses-v-lem} on $v'\sigma$. We conclude using property \eqref{queue-sous-gauss} of the random variable $U_1$.
\end{proof}

\begin{proof} 
By Lemma \ref{lemme-prelim} we prove easily that for every $n \ge 0$
\begin{equation} \label{majo-proba-vois-delta}
	\probc[1]{\Delta \in (X_n, X_{n+1})}{\F_n} \le \exp \pa{-\frac{1}{c_\sigma^2 \gamma_{n+1}}} \quad \text{on ${X_n \in \mathcal{U}_\Delta}$}. 
\end{equation}
We now consider the event $\ac{X_n \notin \mathcal{U}_\Delta}$. Then we have 
\begin{align*}
	\ac{\Delta \in (X_n, X_{n+1})} & = \ac{\exists t \in [0,1], \; U_{n+1} = \frac{\Delta - X_n}{t \sqrt{\gamma_{n+1}} \sigma(X_n)} - \sqrt{\gamma_{n+1}} \frac{b(X_n)}{\sigma(X_n)}}, \\
	& \subset \ac{ \abs{U_{n+1}} \ge \frac{\abs{\Delta - X_n}}{\sqrt{\gamma_{n+1}} \abs{\sigma(X_n)}} - \sqrt{\gamma_{n+1}} \frac{\abs{b(X_n)}}{\abs{\sigma(X_n)}}}.
\end{align*}
As the drift $b$ is dominated by $C_b (1+\abs{x})$, we have
\begin{equation*}
	\ac{\Delta \in (X_n, X_{n+1})} \subset \ac{ \abs{U_{n+1}} \ge  \pa{ \frac{\abs{\Delta - X_n}}{\sqrt{\gamma_{n+1}} C_b (1+\abs{X_n})} - \sqrt{\gamma_{n+1}}} \frac{C_b (1+\abs{X_n})}{\abs{\sigma(X_n)}}},
\end{equation*}
and using the triangular inequality and $\abs{X_n-\Delta} \ge \varepsilon$ we prove that 
$\frac{\abs{\Delta-X_n}}{1+\abs{X_n}} \ge \frac{1}{1+\frac{1+\abs{\Delta}}{\varepsilon}}$.
We also deduce that there exists $n_1 \ge 0$ and $C > 0$ such that for any $n \ge n_1$, 
\begin{align}
	\probc{\ac{\Delta \in (X_n, X_{n+1})} \cap \ac{X_n \notin \mathcal{U}_\Delta}}{\F_n} & \le \probc{\ac{\abs{U_{n+1}} \ge  \frac{C}{\sqrt{\gamma_{n+1}}} \frac{C_b(1+\abs{X_n})}{\abs{\sigma(X_n)}}} \cap\ac{X_n \notin \mathcal{U}_\Delta} }{\F_n}, \notag \\
	& \le \probc{\ac{ \abs{U_{n+1}} \ge \frac{C C_b}{C_\sigma \sqrt{\gamma_{n+1}}}} \cap \ac{X_n \notin \mathcal{U}_\Delta}}{\F_n}, \label{majo-proba-ailleurs}
\end{align}
using $\abs{\sigma} \le \sqrt{C_\sigma} \sqrt{V}$.

From \eqref{majo-proba-vois-delta} and \eqref{majo-proba-ailleurs} combined with \eqref{queue-sous-gauss}, we get
\begin{equation*}
	\exists n_1 \ge 0, \; \exists C > 0, \; \forall n \ge n_1, \quad \probc[1]{\Delta \in (X_n, X_{n+1})}{\F_n} \le \exp \pa{-\frac{C}{\gamma_{n+1}}}.
\end{equation*}
By the condition on the sequence $\seq[1]{\gamma}$ and the conditional Borel-Cantelli Lemma we deduce that the event $\ac{\Delta \in (X_n, X_{n+1})}$ occurs a finite number of times.
\end{proof}

\begin{rmrk} The condition on the step sequence $\seq[1]{\gamma}$ is not restrictive. Indeed, it is satisfied for $\seq[1]{\gamma}$ defined by $\gamma_n = \gamma_0 n^{-r}$ with $\gamma_0 > 0$ and $r \in ]0, 1]$, or $\gamma_n = \log(n)^{-r}$ with $r > 1$.
\end{rmrk}

The technical assumption ``$v$ convex'' is not very restrictive in practice. The important point to note is the condition $v'b \ge 0$ which implies that $\Delta$ is unstable for the ODE $u' = b(u)$. But $\Delta$ may be repulsive as well as attractive for the SDE (cf. Remark \ref{lien_ODE}). The condition $v'\sigma = \domin{v}$ in a neighborhood of $\Delta$ is very important and it seems difficult to relax it.

\begin{rmrk} It is easy to extend the above Theorem to the case of finitely many boundary points $\Delta_i$. If for every $\Delta_i$, there exists a neighborhood $\mathcal{U}_{\Delta_i}$ and a convex function $v_i \in \C^2(\mathcal{U}_{\Delta_i}, \R)$ such that $v(\Delta) = 0$ and satisfying \eqref{hypotheses-v}, then  
\begin{equation*}
	\prob{\exists i \in \ac{0,\dots,l}, \; \exists n_0 \ge 0, \; \forall n \ge n_0, \; X_n \in ]\Delta_i, \Delta_{i+1}[} = 1,
\end{equation*}
with $\Delta_0 = -\infty$ and $\Delta_{l+1} = +\infty$.
\end{rmrk}

\subsection{Weighted empirical measures}

Let $\seq[1]{\eta}$ a positive sequence, called weight sequence, such that $H_n = \sum_{k=1}^n \eta_k$ increases to $+\infty$ when $n$ tends to $+\infty$. We define the weighted empirical measures $\seq[1]{\nu^\eta}$ by 
\begin{equation*}
	\forall n \ge 1, \quad \nu^\eta_n(\d x) = \frac{1}{H_k} \sum_{k=1}^n \eta_k \delta_{X_{k-1}}.
\end{equation*}
In this section we assume that the diffusion satisfies a stability condition \ie 
\begin{equation*}
	\exists \alpha > 0, \; \forall \abs{x} \ge M, \quad x b(x) + \frac{1}{2} \sigma^2(x) \le - \alpha x^2.
\end{equation*}
This condition implies that the points $-\infty$ and $+\infty$ are strongly repulsive and that the empirical measures $\proc{\nu}$ are tight. Since $b$ and $\sigma$ have sublinear growth, this condition implies also the tightness of the weighted empirical measures $\seq[1]{\nu^\eta}$ of the scheme and that any weak limit is an invariant probability for the diffusion (cf. \cite{lambpag-bernoulli} or \cite{these}).

A consequence of Theorem \ref{dim1-thm-euler} is the following proposition which describe the convergence of $\seq[1]{\nu^\eta}$ according to the behavior of $b$ and $\sigma$ in a neighborhood of $\Delta$. For more clearness, we parametrize $b$ and $\sigma$. 

\begin{prpstn} Let $\Delta$ the unique point of $\R$ such that $b(\Delta) = \sigma(\Delta) = 0$. We assume that in a neighborhood of $\Delta$ we have $b(x) = \sgn(x-\Delta) \rho_b(x)$ and $\sigma(x) = \rho_\sigma(x)$ with $\rho_b \ge 0$,  
\begin{equation}
	\rho_b(x) \sim c_b \abs{x-\Delta}^{\beta} \quad \text{and} \quad \sigma(x) \sim c_\sigma \abs{x-\Delta}^\varsigma,
\end{equation}
where $\beta$, $\varsigma$, $c_b$ and $c_\sigma$ are positive real numbers and $\varsigma \ge 1$. If the step sequence $\seq[1]{\gamma}$ satisfies $\forall C > 0$, $\sum_{n \ge 1} \exp(-C / \gamma_n) < +\infty$, then 
\begin{itemize}
\item if $1+\beta-2\varsigma > 0$ then $\Delta$ is an attractive boundary point and $\nu^\eta_n \convet \delta_\Delta$,
\item if $1+\beta-2\varsigma = 0$ and $c_\sigma > \sqrt{2 c_b}$ then $\Delta$ is an attractive boundary point and $\nu^\eta_n \convet \delta_\Delta$,
\item if $1+\beta-2\varsigma = 0$, $c_\sigma < \sqrt{2 c_b}$ and $\beta = 1$ (which implies $\varsigma = 1$) then $\Delta$ is a strongly repulsive boundary point and $\nu^\eta_n \convet \nu_+$ or $\nu^\eta_n \convet \nu_-$,
\item if $1+\beta-2\varsigma < 0$, $c_\sigma \le \sqrt{2 c_b}$ and $\beta \in ]0, 1]$ then $\Delta$ is a strongly repulsive boundary point and $\nu^\eta_n \convet \nu_+$ or $\nu^\eta_n \convet \nu_-$,
\end{itemize}
where $\nu_+$ is the invariant probability on $]\Delta, +\infty[$ and $\nu_-$ the invariant probability on $]-\infty, \Delta[$. 
\end{prpstn}

\begin{proof} To simplify notation, we assume without loss of generality that $\Delta = 0$.

Let $\mathcal{U}_\Delta = ]\Delta-\varepsilon, \Delta+\varepsilon[$ a neighborhood of $\Delta$ and the convex function $v(x) = x^2$. Then for every $x \in \mathcal{U}_\Delta$ we have
\begin{equation*}
	v'b = 2 x \sgn(x) \rho_b(x) \ge 0.	
\end{equation*}
Moreover $v' \sigma \sim c_\sigma \abs{x}^{\varsigma+1}$ with $\varsigma \ge 1$ hence there exists $C > 0$ such that $\abs{v' \sigma} \le C v$. By Theorem \ref{dim1-thm-euler} we know then that the scheme lives in $]-\infty, \Delta[$ or in $]\Delta, +\infty[$ after an almost-surely finite random time.  

Furthermore, with $v(x) = x^2$ we have  
\begin{equation}
	\forall x \in \mathcal{U}_\Delta, \quad \A v(x) = 2 \abs{x} \rho_b(x) + \sigma^2(x) = 2 c_b \abs{x}^{1+\beta} + c_\sigma^2 \abs{x}^{2 \varsigma} + \neglig{\abs{x}^{(1+\beta) \vee (2 \varsigma)}},
\end{equation}
and $\ds \frac{1}{2} \sigma^2(x) \frac{(v'(x))^2}{v(x)} = 2 \sigma^2(x) \sim 2 c_\sigma^2 \abs{x}^{2 \varsigma}$. 

-- If $1+\beta > 2 \varsigma$, there exists a neighborhood of $0$ in which $\A v < \frac{1}{2} \sigma \frac{(v')^2}{v}$. Hence by Corollary \ref{dim1-cor-caract}, $0$ is an attractive boundary point. The sequence of weighted empirical measures $\seq[1]{\nu^\eta}$ of the scheme is tight on $[0, +\infty[$ and on $]-\infty, 0]$, and any weak limit is an invariant probability. However, $\delta_0$ (the Dirac at $0$) is the unique invariant probability on $[0, +\infty[$ or on $]-\infty, 0]$. Thus any weak limit of $\seq[1]{\nu^\eta}$ is $\delta_0$, which proves the first item.

-- If $1+\beta = 2 \varsigma$, $\A v(x) = (2 c_b + c_\sigma^2) \abs{x}^{2 \varsigma} + \neglig{\abs{x}^{2 \varsigma}}$. If $c_\sigma > \sqrt{2 c_b}$ there exists a neighborhood of $0$ in which $\A v < 2 \sigma^2$ and we conclude as above.

If $\beta = 1$ and $c_b < \sqrt{2 c_b}$, we have for every $x$ in $\mathcal{U}_{\Delta}$, $\A v(x) \ge 2 \sigma^2(x) + \varepsilon x^2$. By Corollary \ref{dim1-cor-caract}, the point $0$ is then strongly repulsive. Any weak limit of $\seq{\nu^\eta}$ is a probability on $]-\infty, 0[$ or on $]0, +\infty[$, therefore we have $\nu^\eta_n \convet \nu_-$ or $\nu^\eta_n \convet \nu_+$ (we recall that the boundary points $-\infty$ and $+\infty$ are strongly repulsive).

-- The proof for the case $1+\beta < 2 \varsigma$, $c_\varsigma \le \sqrt{2 c_b}$ and $\beta \in ]0, 1]$ is similar. 
\end{proof}

To illustrate this result we take the same example as in the previous section (Example \ref{exemple-}). 
\begin{xmpl} Let $b$ and $\sigma$ be defined by 
\begin{equation*}
	b(x) = \begin{cases} 
	    - 2 \pa[1]{x - 3 \sgn(x)}& \text{if $\abs{x} \ge 3$}, \\ 
	    - \frac{1}{18} x^3 + \frac{1}{2} x & \text{if $\abs{x} \le 3$}
	    \end{cases},
	\quad \text{and} \quad \sigma(x) = c x \; \text{with $c \in ]0,2[$}.
\end{equation*}
We recall that the nature of the boundary point $0$ for the diffusion $\proc{X}$ depends on parameter $c$. If $c \in ]1,2[$, $0$ is an attractive boundary point, if $c = 1$ it is repulsive and if $c \in ]0, 1[$ it is strongly repulsive.

By the above Proposition, we know that the weighted empirical measure $\seq[1]{\nu^\eta}$ weakly converge to $\delta_0$ when $c \in ]1, 2[$ and to $\nu_+$ or $\nu_-$ when $c \in ]0, 1[$. Note that the convergence to $\nu_+$ or $\nu_-$ does not depend on the initial condition and is not previsible.

We give a representation of the density of $\nu$ approximated by $\nu^\eta_n$ with $n = 10^6$. More precisely, we have discretized the interval $[-2,8]$ using $200$ intervals $I_i$ of length $0.05$ and we have computed $\nu^\eta_n(\un_{I_i})$ for each $I_i$ with $n = 10^6$. The step sequence $\seq[1]{\gamma}$ is defined by $\gamma_n = n^{-1/3}$ and the weight sequence $\seq[1]{\eta}$ is defined by $\eta_n = 1$. The results of this approximation of the stationary density are given in Figure \ref{dim1-figure-piege} for different values of $c$.

\begin{figure}[ht!] 
    \begin{minipage}[c]{0.46\linewidth}
	\subfigure[$c = 0.1$]{
	    \input{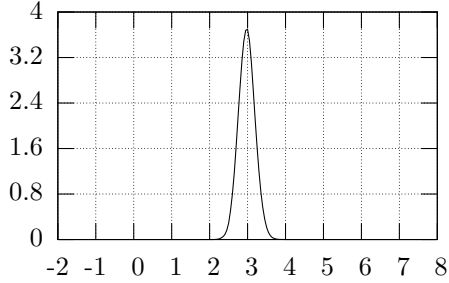}
	}
    \end{minipage}
    \hfill
    \begin{minipage}[c]{0.46\linewidth}
	\subfigure[$c = 0.5$]{ 
	    \input{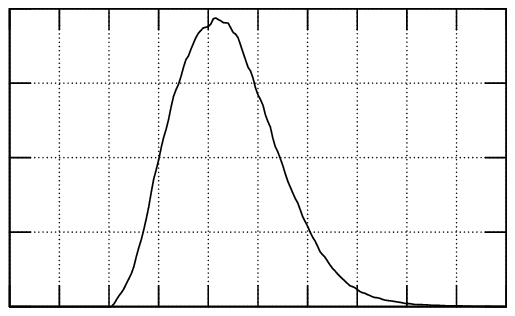}
	}
    \end{minipage}
    
    \begin{minipage}[c]{0.46\linewidth}
	\subfigure[$c = 0.75$]{
	    \input{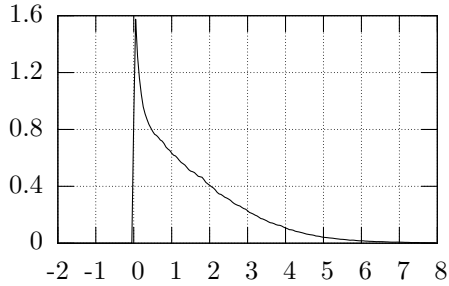}
	}
    \end{minipage}
    \hfill
    \begin{minipage}[c]{0.46\linewidth}
	\subfigure[$c = 1$]{
	    \input{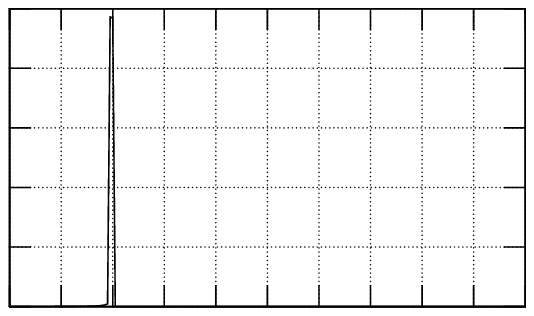}
	}
    \end{minipage}

    \caption{\label{dim1-figure-piege} Approximation of the stationary density for different values of $c$.}
\end{figure}

We remark that for a small noise ($c = 0.1$), the invariant probability concentrates around a stable point of the ODE $u'=b(u)$ (the point $3$), and the more coefficient of diffusion increases, the more the invariant measure is spread out. For $c = 0.75$ we show that the invariant measure is infinite at $0$, and for $c = 1$ the invariant measures seems to be the Dirac mass at $0$.
\end{xmpl}

\section{Numerical example in dimension 2}
To conclude this paper we give a numerical example in dimension 2. We represent the empirical measures of the Euler scheme in a degenerate situation where there are two invariant measures. The first one is the Dirac mass at $\Delta = (0,0)$ and the second one is a probability measure on $\R^2\backslash \ac{(0,0)}$. 

We consider the deterministic Van der Pol equation defined by 
\begin{equation}
	\begin{cases} x' = y, \\
	y' = (1-x^2) y - x.
	\end{cases}
\end{equation}
This non-linear system of $\R^2$ has a stable point $(0,0)$ and an attractive limit cycle. If we add sufficient strong noise in whole space but not in $(0,0)$, the point $(0,0)$ becomes an attractive point for the stochastic system and the invariant measure of this system is $\delta_{(0,0)}$.

Write $b(x, y) =  \begin{pmatrix} y \\ (1-x^2) y - x \end{pmatrix}$ and $\sigma(x, y) = \begin{pmatrix} c x & 0 \\ 0 & c y \end{pmatrix}$. We consider the following perturbed Van der Pol equation 
\begin{equation*}
	\d u_t = b(u_t) \d t + \sigma(u_t) \d B_t
\end{equation*}
with $u_t = (X_t, Y_t) \in \R^2$. We discretize the solution $\proc{X}$ using a decreasing step Euler scheme $\seq{X}$ and the step sequence $\seq[1]{\gamma}$ is defined by $\gamma_n = 0.5 n^{-1/3}$. To guarantee stability of this scheme, we replace the function $b$ by the function $\tilde{b}(x,y) = \begin{pmatrix} y \\ (1-x^2 \wedge 4) y - x \end{pmatrix}$. The scheme is thus defined by $X_0 = (1, 1)$ and for every $n \ge 0$ 
\begin{equation*}
	X_{n+1} = X_n + \gamma_n \tilde{b}(X_n) + \sqrt{\gamma_n} \sigma(X_n) U_{n+1}, 
\end{equation*}
where $U_n$ is a normalized Gaussian of $\R^2$.

\begin{centering}
\begin{figure}[ht!]
	\begin{minipage}[c]{0.8\linewidth}
	\input{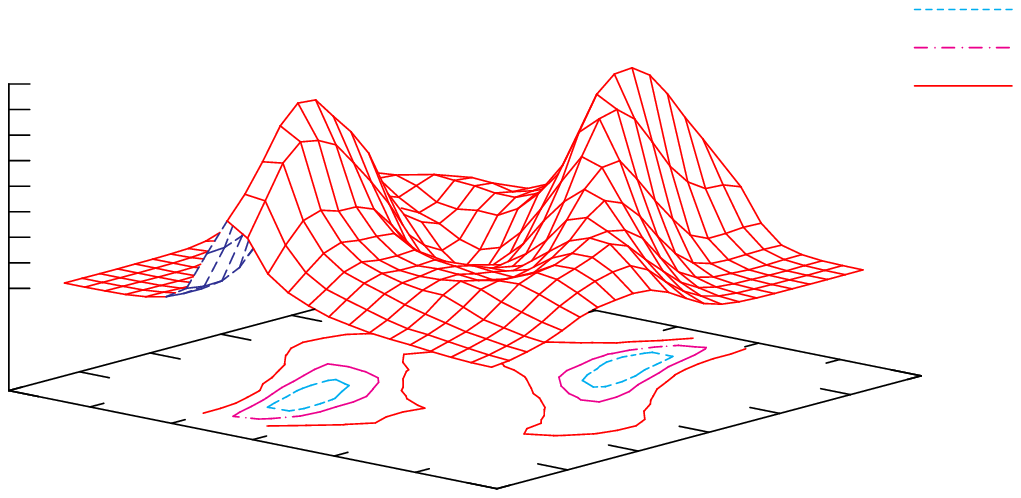}
	\end{minipage}
	\hfill (a) $c = 0.5$
	
	\vspace{-1.7cm}
	\begin{minipage}[c]{0.8\linewidth}
	\input{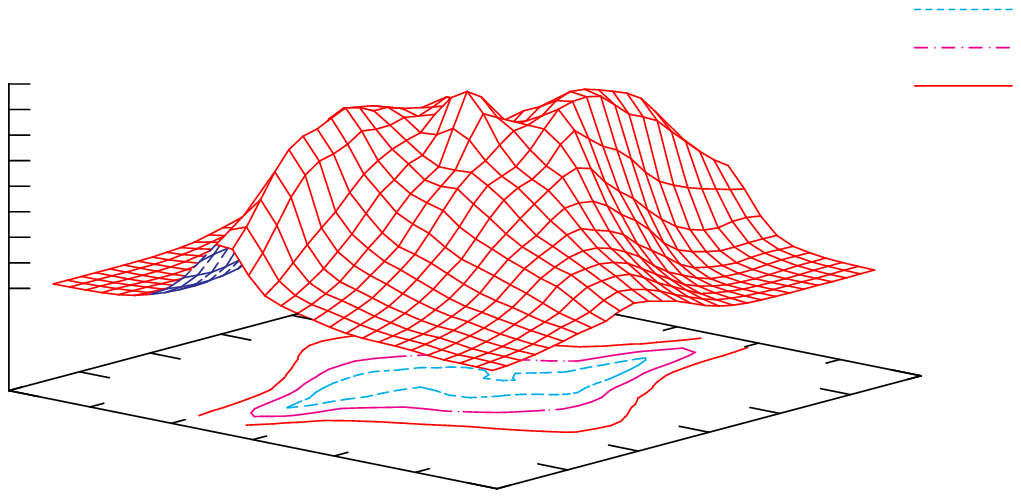}
	\end{minipage}
	\hfill (b) $c = 0.7$
	
	\vspace{-1.7cm}
	\begin{minipage}[c]{0.8\linewidth}
	\input{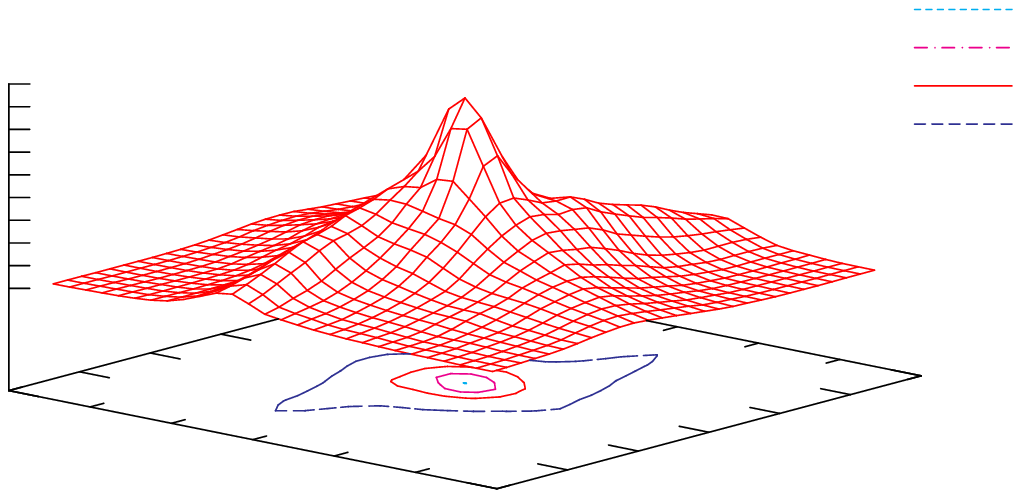}
	\end{minipage}
	\hfill (c) $c = 0.8$
	\caption{\label{dim1-figure-vdp} Approximation of the stationary density of the perturbed Van der Pol equation for different values of $c$.}   
\end{figure}
\end{centering}

An approximation of the density of $\nu$ is done using an histogram of $\nu^\eta_n$ for $n = 10^6$ and $\eta_n = 1$. The histogram is built using a step $h = 0.2$. The results are presented in Figure \ref{dim1-figure-vdp}.

We note that for a small noise ($c = 0.5$) the point $(0,0)$ seems not charged by the invariant probability. And for higher values of $c$ the invariant probability concentrates in a neighborhood of $(0,0)$.

\bibliography{biblio-article,biblio-livre}
\bibliographystyle{plain}

\end{document}